\documentclass[final,3p]{elsarticle}
 \usepackage{graphics}
 \usepackage{graphicx}
 \usepackage{epsfig}
\usepackage{amssymb}
 \usepackage{amsthm}
 \usepackage{lineno}
 \usepackage{amsmath}
   \numberwithin{equation}{section}
\usepackage{mathrsfs}

\newtheorem{thm}{Theorem}[section]
\newtheorem{cor}[thm]{Corollary}
\newtheorem{lem}[thm]{Lemma}
\newtheorem{prop}[thm]{Proposition}
\newtheorem{defn}[thm]{Definition}
\newtheorem{rem}[thm]{Remark}

\biboptions{sort&compress,square}
\begin{document}
\begin{frontmatter}
 \author{Jian Wang}
 \ead{wangj068@gmail.com}
\author{Yong Wang\corref{cor2}}
\cortext[cor2]{Corresponding author;  ~Email address: \textbf{wangy581@nenu.edu.cn (Yong Wang)}}
\address{School of Mathematics and Statistics, Northeast Normal University,
Changchun, 130024, P.R.China}
\title{The Kastler-Kalau-Walze type
 \\theorem for 6-dimensional manifolds with boundary}
\begin{abstract}
In this paper, we define lower dimensional
volumes of spin manifolds with boundary. We compute the lower
dimensional volume ${\rm Vol}_{6}^{(1,3)}$ for $6$-dimensional spin manifolds with boundary
and the gravity on boundary is derived by the noncommutative residue associated with Dirac operators.
For 6-dimensional manifolds with boundary, we also get a Kastler-Kalau-Walze type theorem for a general fourth order operator .
\end{abstract}
\end{frontmatter}
\section{Introduction}
The noncommutative residue plays a prominent role in noncommutative geometry \cite{Gu}\cite{Wo}. Connes \cite{Co1} used the noncommutative
 residue to derive a conformal 4-dimensional Polyakov action analogy. Connes \cite{Co2} proved that the noncommutative
residue on a compact manifold $M$ coincided with the Dixmier's trace on pseudodifferential operators of order $-{\rm {dim}}M$. Several
years ago, Connes made a challenging observation that the noncommutative residue of the square of the inverse of the Dirac
operator was proportional to the Einstein-Hilbert action, which we call the Kastler-Kalau-Walze theorem. Kastler\cite{Ka} gave a
brute-force proof of this theorem. Kalau and Walze \cite{KW} proved this theorem in the normal coordinates system simultaneously.
Ackermann \cite{Ac} gave a note on a new proof of this theorem by means of the heat kernel expansion.

  Recently, Ponge defined lower dimensional volumes of Riemannian manifolds by the Wodzicki residue \cite{Po}.
Fedosov et al defined a noncommutative residue on Boutet de Monvel's algebra and proved that it was a
unique continuous trace \cite{FGLS}. Wang generalized  the Connes' results to the case of manifolds with boundary in \cite{Wa1} \cite{Wa2} ,
and proved a Kastler-Kalau-Walze type theorem for the Dirac operator and the signature operator for lower-dimensional manifolds
with boundary\cite{Wa3} \cite{Wa4}. The purpose of papers \cite{Wa3} \cite{Wa4} is to derive the gravitational action by the noncommutative
 residue  associated with  Dirac operators for spin manifolds with boundary, but the boundary term vanished in Wang's results.
The motivation of this paper is to derive the gravitational action on boundary by the noncommutative residue  associated with
 Dirac operators. In other words, we want to get a nonvanishing boundary term. In \cite{Wa3} \cite{Wa4}, Wang computed
 $\widetilde{{\rm Wres}}[\pi^+D^{-1}\circ\pi^+D^{-1}]$ and $\widetilde{{\rm Wres}}[\pi^+D^{-2}\circ\pi^+D^{-2}]$, where the two
 operators are symmetric. In this paper,  for 6-dimensional manifolds with boundary, we compute
  $\widetilde{{\rm Wres}}[\pi^+D^{-1}\circ\pi^+D^{-3}]$. Since $D^{-1}$ and $D^{-3}$ are not symmetric, the
  boundary term is the integral of the extrinsic scalar curvature and the gravitational action on boundary emerges.
 We also get a generalized Kastler-Kalau-Walze theorem associated with a general fourth order operators  for
 $6$-dimensional spin manifolds.

 This paper is organized as follows: In Section 2, we define lower dimensional volumes of spin manifolds with boundary. In
Section 3, for $6$-dimensional spin manifolds with boundary and the associated Dirac operator $D$ and $D^{3}$,
 we compute the lower dimensional volume ${\rm Vol}^{(1,3)}_6$ and get a Kastler-Kalau-Walze type theorem in this case.
In Section 4, we get a Kastler-Kalau-Walze type theorem associated with the fourth order operators  for
 $6$-dimensional spin manifolds.

\section{Lower dimensional volumes of spin manifolds with boundary}
 In order to define lower dimensional volumes of spin manifolds with boundary, we need some basic facts and formulae about Boutet de
Monvel's calculus and the definition of the noncommutative residue for manifolds with boundary. We can find them in Section 2,3 \cite{Wa2}
and Section 2.1\cite{Wa3}.

Let $M$ be a $n$-dimensional compact oriented spin manifold with boundary $\partial M$. We assume that the metric $g^M$ on $M$
has the following form near the boundary,
 \begin{equation}
 g^M=\frac{1}{h(x_n)}g^{\partial M}+dx_n^2,
\end{equation}
where $g^{\partial M}$ is the metric on  ${\partial M}$. Let $D$ be the Dirac operator associated to $g$ on the spinors bundle $S(TM)$\cite{Wa3}.
 Let $p_1,p_2$ be nonnegative integers and $p_1+p_2\leq n$.

\begin{defn}
Lower dimensional volumes of spin manifolds with boundary are defined by
 \begin{equation}
{\rm Vol}^{(p_1,p_2)}_nM:= \widetilde{{\rm Wres}}[\pi^+D^{-p_1}\circ\pi^+D^{-p_2}]
\end{equation}
where the related definitions, see Section 2, 3\cite{Wa3}.
\end{defn}
 Denote by $\sigma_l(A)$ the $l$-order symbol of an operator $A$. By (2.1.4)-(2.1.8)\cite{Wa3}, we get
\begin{equation}
\widetilde{{\rm Wres}}[\pi^+D^{-p_1}\circ\pi^+D^{-p_2}]=\int_M\int_{|\xi|=1}{\rm
trace}_{S(TM)}[\sigma_{-n}(D^{-p_1-p_2})]\sigma(\xi)dx+\int_{\partial M}\Phi,
\end{equation}
and
\begin{eqnarray}
\Phi &=&\int_{|\xi'|=1}\int^{+\infty}_{-\infty}\sum^{\infty}_{j, k=0}\sum\frac{(-i)^{|\alpha|+j+k+1}}{\alpha!(j+k+1)!}
\times {\rm trace}_{S(TM)}[\partial^j_{x_n}\partial^\alpha_{\xi'}\partial^k_{\xi_n}\sigma^+_{r}(D^{-p_1})(x',0,\xi',\xi_n)
\nonumber\\
&&\times\partial^\alpha_{x'}\partial^{j+1}_{\xi_n}\partial^k_{x_n}\sigma_{l}(D^{-p_2})(x',0,\xi',\xi_n)]d\xi_n\sigma(\xi')dx',
\end{eqnarray}
 where the sum is taken over $r-k-|\alpha|+l-j-1=-n,~~r\leq -p_1,l\leq -p_2$.

 Since $[\sigma_{-n}(D^{-p_1-p_2})]|_M$ has the same expression as $\sigma_{-n}(D^{-p_1-p_2})$ in the case of manifolds without
boundary, so locally we can use the computations \cite{Ka}, \cite{KW}, \cite{Po}, \cite{Wa3} to compute the first term.
The following proposition is the motivation
of the definition of lower dimensional volumes of spin manifolds with boundary \cite{Wa4}.

\begin{prop}\cite{Wa4}  Lower dimensional volumes of spin manifolds with boundary are given by

1) When $p_1+p_2=n$, then ${\rm Vol}^{(p_1,p_2)}_nM=c_0{\rm Vol}_M.$

2) when $p_1+p_2\equiv n {\rm mod} 1$, ${\rm Vol}^{(p_1,p_2)}_nM=\int_{\partial M}\Phi.$

3)  \begin{equation}
{\rm Vol}^{(1,1)}_4=-\frac{\Omega_4}{3}\int_Ms{\rm dvol}_M;~~ {\rm Vol}^{(1,1)}_3=c_1{\rm Vol}_{\partial M}
\end{equation}
 where $c_0,c_1$ are constants and $s$ is the scalar curvature.
\end{prop}

\section{A Kastler-Kalau-Walze type theorem for $6$-dimensional spin manifolds with boundary }
 In this section, We compute the lower dimensional volume ${\rm Vol}^{(1,3)}_6$ for $6$-dimensional spin manifolds with
boundary and get a Kastler-Kalau-Walze type theorem in this case.

Firstly, we compute the symbol $\sigma(D^{-3})$ of $D^{-3}$. Recall the definition of the Dirac operator $D$  \cite{FGLS}\cite{Y}.
 Let $\nabla^L$ denote the Levi-Civita connection about $g^M$. In the local coordinates $\{x_i; 1\leq i\leq n\}$ and the
 fixed orthonormal frame $\{\widetilde{e_1},\cdots,\widetilde{e_n}\}$, the connection matrix $(\omega_{s,t})$ is defined by
\begin{equation}
\nabla^L(\widetilde{e_1},\cdots,\widetilde{e_n})= (\widetilde{e_1},\cdots,\widetilde{e_n})(\omega_{s,t}).
\end{equation}

 The Dirac operator is defined by
\begin{equation}
D=\sum^n_{i=1}c(\widetilde{e_i})\Big[\widetilde{e_i}
-\frac{1}{4}\sum_{s,t}\omega_{s,t}(\widetilde{e_i})c(\widetilde{e_s})c(\widetilde{e_t})\Big].
\end{equation}
where $c(\widetilde{e_i})$ denotes the Clifford action.

Recall the definition of the Dirac operator $D^{2}$ in \cite{Ka}, \cite{KW} and \cite{Wa4}, we have
\begin{equation}
D^{2}=-\sum_{i,j}g^{i,j}\Big[\partial_{i}\partial_{j}+2\sigma_{i}\partial_{j}+(\partial_{i}\sigma_{j})+\sigma_{i}\sigma_{j}
    -\Gamma_{i,j}^{k}\partial_{i}-\Gamma_{i,j}^{k}\sigma_{k}\Big]+\frac{1}{4}s.
\end{equation}
where $\sigma_{i}:=-\frac{1}{4}\sum_{s,t}\omega_{s,t}(\partial_{i})e_se_t$.

Combining (3.2) and (3.3), we have
\begin{eqnarray}
D^{3}
&=&\sum^n_{i=1}c(\widetilde{e_l}) \langle e_l, dx_{l}\rangle \bigg\{ (-\sum_{i,j}g^{i,j}\partial_{i}\partial_{j})\partial_{l}
 -\sum_{i,j}\partial_{l} (g^{i,j})\partial_{i}\partial_{j} \nonumber\\
&& +
 \Big[-\sum_{i,j}g^{i,j}(2\sigma_{i}\partial_{j}-\Gamma_{i,j}^{k}\partial_{i})\Big]\partial_{l}
-\sum_{i,j}\partial_{l}(g^{i,j})(2\sigma_{i}\partial_{j}-\Gamma_{i,j}^{k}\partial_{i})   \nonumber\\
&&-\sum_{i,j}g^{i,j}\Big(2(\partial_{l}\sigma_{i})\partial_{j}
 -\partial_{l}(\Gamma_{i,j}^{k})\partial_{i}\Big)  \nonumber\\
&& +
  \Big[-\sum_{i,j}g^{i,j}\Big((\partial_{i}\sigma_{j})+\sigma_{i}\sigma_{j}-\Gamma_{i,j}^{k}\sigma_{k}\Big)+\frac{1}{4}s\Big]\partial_{l}   \nonumber\\
&& -\sum_{i,j}\partial_{l}(g^{i,j})\Big((\partial_{i}\sigma_{j})+\sigma_{i}\sigma_{j}
  -\Gamma_{i,j}^{k}\sigma_{k}\Big)+\frac{1}{4}\partial_{l}s
    \nonumber\\
&& -\sum_{i,j}g^{i,j}\partial_{l}\Big((\partial_{i}\sigma_{j})
+\sigma_{i}\sigma_{j}-\Gamma_{i,j}^{k}\sigma_{k}\Big)+\frac{1}{4}\partial_{l}s\bigg\}
     \nonumber\\
&& -\frac{1}{4}\sum_{s,t}\omega_{s,t}(\widetilde{e_l})c(e_{l})c(\widetilde{e_s})c(\widetilde{e_t})
   \bigg\{-\sum_{i,j}g^{i,j}\partial_{i}\partial_{j}
   -\sum_{i,j}g^{i,j}(2\sigma_{i}\partial_{j}-\Gamma_{i,j}^{k}\partial_{i})  \nonumber\\
&& -\sum_{i,j}g^{i,j}\Big((\partial_{i}\sigma_{j})+\sigma_{i}\sigma_{j}-\Gamma_{i,j}^{k}\sigma_{k}\Big)+\frac{1}{4}s\bigg\}.
\end{eqnarray}
Then, we obtain
\begin{eqnarray}
\sigma_{3}(D^{3})&=&\sqrt{-1}c(\xi)|\xi|^2 , \\
\sigma_{2}(D^{3})&=&c(\xi)(2\sigma^k-\Gamma^k)\xi_{k}
     -\frac{1}{4}|\xi|^2\sum_{s,t}\omega_{s,t}(\widetilde{e_l})c(e_{l})c(\widetilde{e_s})c(\widetilde{e_t}), \\
\sigma_{1}(D^{3})&=&\sum^n_{i=1}c(\widetilde{e_l}) \langle e_l, dx_{l}\rangle
    \bigg\{ -\sum_{i,j}\partial_{l}(g^{i,j})(2\sigma_{i}\partial_{j}-\Gamma_{i,j}^{k}\partial_{i})    \nonumber\\
      &&+\sum_{i,j}g^{i,j}\partial_{l}(\Gamma_{i,j}^{k})\partial_{i}
     -\sum_{i,j}g^{i,j}\Big(2(\partial_{l}\sigma_{i})\partial_{j}\Big)   \nonumber\\
 && + \Big[-\sum_{i,j}g^{i,j}\Big((\partial_{i}\sigma_{j})+\sigma_{i}\sigma_{j}
 -\Gamma_{i,j}^{k}\sigma_{k}\Big)+\frac{1}{4}s\Big]\partial_{l}\bigg\}
     \nonumber\\
      && -\frac{1}{4}\sum_{s,t}\omega_{s,t}(\widetilde{e_l})c(e_{l})c(\widetilde{e_s})c(\widetilde{e_t})
           \Big[-\sum_{i,j}g^{i,j}(2\sigma_{i}\partial_{j}-\Gamma_{i,j}^{k}\partial_{i})\Big],   \\
\sigma_{0}(D^{3})&=&
     \sum^n_{i=1}c(\widetilde{e_l}) \langle e_l, dx_{l}\rangle
  \bigg\{-\sum_{i,j}\partial_{l}(g^{i,j})\Big((\partial_{i}\sigma_{j})+\sigma_{i}\sigma_{j}-\Gamma_{i,j}^{k}\sigma_{k}\Big)+\frac{1}{4}
    \partial_{l}s   \nonumber\\
&& -\sum_{i,j}g^{i,j}\partial_{l}\Big((\partial_{i}\sigma_{j})+\sigma_{i}\sigma_{j}-\Gamma_{i,j}^{k}\sigma_{k}\Big)+\frac{1}{4}
   \partial_{l}s \bigg\}    \nonumber\\
&& -\frac{1}{4}\sum_{s,t}\omega_{s,t}(\widetilde{e_l})c(e_{l})c(\widetilde{e_s})c(\widetilde{e_t})
\Big[-\sum_{i,j}g^{i,j}\Big((\partial_{i}\sigma_{j})+\sigma_{i}\sigma_{j}-\Gamma_{i,j}^{k}\sigma_{k}\Big)+\frac{1}{4}s\Big].
\end{eqnarray}

Write
 \begin{equation}
D_x^{\alpha}=(-\sqrt{-1})^{|\alpha|}\partial_x^{\alpha};
~\sigma(D^{3})=p_3+p_2+p_1+p_0;
~\sigma(D^{-3})=\sum^{\infty}_{j=3}q_{-j}.
\end{equation}
 By the composition formula of psudodifferential operators, we have
 \begin{eqnarray}
1=\sigma(D^{3}\circ D^{-3})&=&\sum_{\alpha}\frac{1}{\alpha!}\partial^{\alpha}_{\xi}[\sigma(D)]D^{\alpha}_{x}[\sigma(D^{-1})] \nonumber\\
&=&(p_3+p_2+p_1+p_0)(q_{-3}+q_{-4}+q_{-5}+\cdots) \nonumber\\
&&+\sum_j(\partial_{\xi_j}p_3+\partial_{\xi_j}p_2+\partial_{\xi_j}p_1+\partial_{\xi_j}p_0)
(D_{x_j}q_{-3}+D_{x_j}q_{-4}+D_{x_j}q_{-5}+\cdots) \nonumber\\
&=&p_3q_{-3}+(p_3q_{-4}+p_2q_{-3}+\sum_j\partial_{\xi_j}p_3D_{x_j}q_{-3})+\cdots,
\end{eqnarray}
Then we obtain
\begin{equation}
q_{-3}=p_3^{-1};~q_{-4}=-p_3^{-1}[p_2p_3^{-1}+\sum_j\partial_{\xi_j}p_3D_{x_j}(p_3^{-1})].
\end{equation}
By Lemma 2.1 in \cite{Wa3} and (3.4)-(3.11), we obtain

\begin{lem}
\begin{eqnarray}
\sigma_{-1}(D^{-1})&=&\frac{\sqrt{-1}c(\xi)}{|\xi|^2}; \\
\sigma_{-2}(D^{-1})&=&\frac{c(\xi)\sigma_{0}(D)c(\xi)}{|\xi|^4}+\frac{c(\xi)}{|\xi|^6}\sum_jc(dx_j)
\Big(\partial_{x_j}[c(\xi)]|\xi|^2-c(\xi)\partial_{x_j}(|\xi|^2)\Big) ;\\
\sigma_{-3}(D^{-3})&=&\frac{\sqrt{-1}c(\xi)}{|\xi|^4}; \\
\sigma_{-4}(D^{-3})&=&\frac{c(\xi)\sigma_{2}(D^{3})c(\xi)}{|\xi|^8}+\frac{\sqrt{-1}c(\xi)}{|\xi|^8}\Big(|\xi|^4c(dx_n)\partial_{x_n}c(\xi')
-2h'(0)c(dx_n)c(\xi)\nonumber\\
&&+2\xi_{n}c(\xi)\partial_{x_n}c(\xi')+4\xi_{n}h'(0)\Big),
\end{eqnarray}
where $\sigma_{0}(D)=-\frac{1}{4}\sum_{s,t}\omega_{s,t}(\widetilde{e_i})c(\widetilde{e_i})c(\widetilde{e_s})c(\widetilde{e_t})$.
\end{lem}

 Since $\Phi$ is a global form on $\partial M$, so for any fixed point $x_0\in\partial M$, we can choose the normal coordinates
$U$ of $x_0$ in $\partial M$ (not in $M$) and compute $\Phi(x_0)$ in the coordinates $\widetilde{U}=U\times [0,1)\subset M$ and the
metric $\frac{1}{h(x_n)}g^{\partial M}+dx_n^2.$ The dual metric of $g^M$ on $\widetilde{U}$ is ${h(x_n)}g^{\partial M}+dx_n^2.$  Write
$g^M_{ij}=g^M(\frac{\partial}{\partial x_i},\frac{\partial}{\partial x_j});~ g_M^{ij}=g^M(dx_i,dx_j)$, then

\begin{equation}
[g^M_{i,j}]= \left[\begin{array}{lcr}
  \frac{1}{h(x_n)}[g_{i,j}^{\partial M}]  & 0  \\
   0  &  1
\end{array}\right];~~~
[g_M^{i,j}]= \left[\begin{array}{lcr}
  h(x_n)[g^{i,j}_{\partial M}]  & 0  \\
   0  &  1
\end{array}\right],
\end{equation}
and
 \begin{equation}
\partial_{x_s}g_{ij}^{\partial M}(x_0)=0, 1\leq i,j\leq n-1; ~~~g_{ij}^M(x_0)=\delta_{ij}.
\end{equation}

Let $n=6$ and $\{e_1,\cdots,e_{n-1}\}$ be an orthonormal frame field in $U$ about $g^{\partial M}$ which is parallel along geodesics and
$e_i(x_0)=\frac{\partial}{\partial x_i}(x_0)$, then $\{\widetilde{e_1}=\sqrt{h(x_n)}e_1,\cdots,\widetilde{e_{n-1}}=\sqrt{h(x_n)}e_{n-1},
\widetilde{e_n}=dx_n\}$ is the orthonormal frame field in $\widetilde U$ about $g^M$. Locally $S(TM)|_{\widetilde {U}}\cong
\widetilde {U}\times\wedge^* _{\bf C}(\frac{n}{2}).$
Let $\{f_1,\cdots,f_8\}$ be the orthonormal basis of $\wedge^* _{\bf C}(\frac{n}{2})$.
Take a spin frame field $\sigma:~\widetilde {U}\rightarrow {\rm Spin}(M)$ such that $\pi\sigma=
\{\widetilde{e_1},\cdots,\widetilde{e_n}\}$, where $\pi :~{\rm Spin}(M)\rightarrow O(M)$ is a double covering,
 then $\{[(\sigma,f_i)],~1\leq i\leq 8\}$ is an orthonormal frame of $S(TM)|_{\widetilde {U}}.$ In the following, since the global form
$\Phi$ is independent of the choice of the local frame, so we can compute ${\rm tr}_{S(TM)}$ in the frame $\{[(\sigma,f_i)],~1\leq
i\leq 8\}.$ Let $\{E_1,\cdots,E_n\}$ be the canonical basis of
${\bf R}^n$ and $c(E_i)\in {\rm cl}_{\bf C}(n)\cong {\rm Hom}(\wedge^*_{\bf C}(\frac{n}{2}),\wedge^* _{\bf C}(\frac{n}{2}))$
be the Clifford action. By [Y], then
\begin{equation}
c( \widetilde{e_i})=[(\sigma,c(E_i))];~ c( \widetilde{e_i})[(\sigma,f_i)]=[(\sigma,c(E_i)f_i)];~
\frac{\partial}{\partial x_i}=[(\sigma,\frac{\partial}{\partial
x_i})],
\end{equation}
then we have $\frac{\partial}{\partial x_i}c( \widetilde{e_i})=0$ in the above frame.
By Lemma 2.2 in \cite{Wa3}, we have
\begin{lem}\cite{Wa3}  For $n$-dimensional spin manifolds with boundary,
\begin{eqnarray}
&&\partial_{x_j}(|\xi|_{g^M}^2)(x_0)=0,~{\rm if}~j<n;~=h'(0)|\xi'|_{g^{\partial M}}^2,~{\rm if }~j=n. \\
&&\partial_{x_j}[c(\xi)](x_0)=0,~{\rm if }~j<n;~=\partial_{x_n}[c(\xi')](x_0),~{\rm if}~j=n,
\end{eqnarray}
 where $\xi=\xi'+\xi_ndx_n$.
\end{lem}

Next we  compute $\sigma_{0}(D)(x_0)$.
By Lemma 2.3 in \cite{Wa3}, we have
\begin{lem}\cite{Wa3}  For $6$-dimensional spin manifolds with boundary,
 \begin{equation}
\sigma_{0}(D)(x_0)=-\frac{5}{4}h'(0)c(dx_n).
\end{equation}
\end{lem}

 Now we can compute $\Phi$ (see formula (2.4) for the definition of $\Phi$), since the sum is taken over $
-r-l+k+j+|\alpha|=5,~~r\leq -1,l\leq-4,$ then we have the following five cases:

{\bf case a)~I)}~$r=-1,~l=-3~k=j=0,~|\alpha|=1$

From (2.4) we have
 \begin{equation}
{\rm case~a)~I)}=-\int_{|\xi'|=1}\int^{+\infty}_{-\infty}\sum_{|\alpha|=1}
 {\rm trace}[\partial^\alpha_{\xi'}\pi^+_{\xi_n}\sigma_{-1}(D^{-1})\times
 \partial^\alpha_{x'}\partial_{\xi_n}\sigma_{-3}(D^{-3})](x_0)d\xi_n\sigma(\xi')dx'.
\end{equation}

By Lemma 3.2, for $i<n$, then
\begin{eqnarray}
\partial_{x_i}\sigma_{-3}(D^{-3})(x_0)&=&\partial_{x_i}[\frac{\sqrt{-1}c(\xi)}{|\xi|^4}](x_0) \nonumber\\
&=&\frac{\sqrt{-1}\partial_{x_i}[c(\xi)](x_0)}{|\xi|^4}-\frac{2\sqrt{-1}c(\xi)\partial_{x_i}[|\xi|^2](x_0)}{|\xi|^6}=0.
\end{eqnarray}
Then case a) I) vanishes.

 {\bf case a)~II)}~$r=-1,~l=-3 ~k=|\alpha|=0,~j=1$

From (2.4) we have
 \begin{equation}
{\rm case~ a)~II)}=-\frac{1}{2}\int_{|\xi'|=1}\int^{+\infty}_{-\infty} {\rm
trace} [\partial_{x_n}\pi^+_{\xi_n}\sigma_{-1}(D^{-1})\times
\partial_{\xi_n}^2\sigma_{-3}(D^{-3})](x_0)d\xi_n\sigma(\xi')dx'.
\end{equation}
By (2.2.23) in \cite{Wa3}, we have
 \begin{equation}
\pi^+_{\xi_n}\partial_{x_n}\sigma_{-1}(D^{-1})(x_0)|_{|\xi'|=1}=\frac{\partial_{x_n}[c(\xi')](x_0)}{2(\xi_n-i)}+\sqrt{-1}h'(0)
\Big[\frac{ic(\xi')}{4(\xi_n-i)}+\frac{c(\xi')+ic(dx_n)}{4(\xi_n-i)^2}\Big].
\end{equation}
By (3.14), we obtain
 \begin{equation}
\partial^2_{\xi_n}\sigma_{-3}(D^{-3})=\sqrt{-1}\Big(\frac{(20\xi_n^{2}-4)c(\xi')+12(\xi_n^{3}-\xi_n)c(dx_n)}{(1+\xi_n^{2})^4}\Big);
\end{equation}
Since $n=6$, ${\rm tr}_{S(TM)}[{\rm id}]={\rm dim}(\wedge^*(3))=8$.
By the relation of the Clifford action and ${\rm tr}{AB}={\rm tr }{BA}$, then we have the equalities:
\begin{eqnarray}
&&{\rm tr}[c(\xi')c(dx_n)]=0;~~{\rm tr}[c(dx_n)^2]=-8;~~{\rm tr}[c(\xi')^2](x_0)|_{|\xi'|=1}=-8; \nonumber\\
&&{\rm tr}[\partial_{x_n}c(\xi')c(dx_n)]=0;~~{\rm tr}[\partial_{x_n}c(\xi')c(\xi')](x_0)|_{|\xi'|=1}=-4h'(0).
\end{eqnarray}
By (3.25)-(3.27) and direct computations ,we have
\begin{eqnarray}
&&{\rm tr}\Big\{\Big[\frac{\partial_{x_n}[c(\xi')](x_0)}{2(\xi_n-i)}+\sqrt{-1}h'(0)
\Big(\frac{ic(\xi')}{4(\xi_n-i)}+\frac{c(\xi')+ic(dx_n)}{4(\xi_n-i)^2}\Big)\Big]\nonumber\\
&&\times
\Big[\sqrt{-1}\Big(\frac{(20\xi_n^{2}-4)c(\xi')+12(\xi_n^{3}-\xi_n)c(dx_n)}{(1+\xi_n^{2})^4}\Big)\Big]
\Big\}(x_0)|_{|\xi'|=1} \nonumber\\
&=&h'(0)\frac{-8-24i \xi_n+40\xi_n^2+24 i\xi_n^{3}}{(\xi_n-i)^6(\xi_n+i)^4}.
\end{eqnarray}
By (3.24) and (3.28), we obtain
\begin{eqnarray}
{\rm case~ a)~II)}&=&\int_{|\xi'|=1}\int^{+\infty}_{-\infty}h'(0)\frac{4+12i\xi_n-20\xi_n^2-12i\xi_n^{3}}{(\xi_n-i)^6(\xi_n+i)^4}
d\xi_n\sigma(\xi')dx' \nonumber\\
&=&h'(0)\Omega_4\int_{\Gamma^+}\frac{4+12i\xi_n-20\xi_n^2-12i\xi_n^{3}}{(\xi_n-i)^6(\xi_n+i)^4}d\xi_ndx'\nonumber\\
&=&h'(0)\Omega_4 \frac{2\pi i}{5!}\Big[\frac{4+12i\xi_n-20\xi_n^2-12i\xi_n^{3}}{(\xi_n+i)^4}\Big]^{(5)}|_{\xi_n=i}dx'\nonumber\\
&=&-\frac{15}{16}\pi h'(0)\Omega_4dx',
\end{eqnarray}
 where $\Omega_4$ is the canonical volume of $S^4$.

{\bf case a)~III)}~$r=-1,~l=-3~j=|\alpha|=0,~k=1$\\

From (2.4) we have
 \begin{equation}
{\rm case~ a)~III)}=-\frac{1}{2}\int_{|\xi'|=1}\int^{+\infty}_{-\infty}
{\rm trace} [\partial_{\xi_n}\pi^+_{\xi_n}\sigma_{-1}(D^{-1})\times
\partial_{\xi_n}\partial_{x_n}\sigma_{-3}(D^{-3})](x_0)d\xi_n\sigma(\xi')dx'.
\end{equation}
By (2.2.29) in \cite{Wa3}, we have
 \begin{equation}
\partial_{\xi_n}\pi^+_{\xi_n}\sigma_{-1}(D^{-1})(x_0)|_{|\xi'|=1}=-\frac{c(\xi')+ic(dx_n)}{2(\xi_n-i)^2}.
\end{equation}
By (3.14), we obtain
 \begin{equation}
\partial_{\xi_n}\partial_{x_n}\sigma_{-3}(D^{-3})(x_0)|_{|\xi'|=1}=-2\sqrt{-1}h'(0)
\Big[\frac{(1-5\xi_n^{2})c(dx_n)}{(1+\xi_n^{2})^4}-\frac{6\xi_n c(\xi')}{(1+\xi_n^{2})^4}\Big]-
\frac{4\sqrt{-1}\xi_n\partial_{x_n}c(\xi')(x_0)}{(1+\xi_n^{2})^3}.
\end{equation}
By (3.27), (3.31) and (3.32), we obtain
\begin{eqnarray}
&&{\rm  tr}\Big\{\Big[\frac{c(\xi')+ic(dx_n)}{2(\xi_n-i)^2}\Big]\times
\Big[2\sqrt{-1}h'(0)
\Big(\frac{(1-5\xi_n^{2})c(dx_n)}{(1+\xi_n^{2})^4}-\frac{6\xi_n c(\xi')}{(1+\xi_n^{2})^4}\Big)
+\frac{4\sqrt{-1}\xi_n\partial_{x_n}c(\xi')(x_0)}{(1+\xi_n^{2})^3}\Big]
\Big\}(x_0)|_{|\xi'|=1} \nonumber\\
&=&h'(0)\frac{8i-32\xi_n-8i\xi_n^{2}}{(\xi_n-i)^5(\xi_n+i)^4}.
\end{eqnarray}
Then
\begin{eqnarray}
{\rm {\bf case~a)~III)}}&=&
 h'(0)\int_{|\xi'|=1}\int^{+\infty}_{-\infty}\frac{-4i+16\xi_n+4i\xi_n^{2}}{(\xi_n-i)^5(\xi_n+i)^4}d\xi_n\sigma(\xi')dx' \nonumber\\
&=& h'(0)\frac{2 \pi i}{4!}\Big[\frac{-4i+16\xi_n+4i\xi_n^{2}}{(\xi_n+i)^4} \Big]^{(4)}|_{\xi_n=i}\Omega_4dx'\nonumber\\
&=&\frac{25}{16}\pi h'(0)\Omega_4dx'.
\end{eqnarray}

{\bf case b)}~$r=-1,~l=-4,~k=j=|\alpha|=0$\\

From (2.4) and an integration by parts, we have
\begin{eqnarray}
{\rm case~ b)}&=&-i\int_{|\xi'|=1}\int^{+\infty}_{-\infty}{\rm trace} [\pi^+_{\xi_n}\sigma_{-1}(D^{-1})\times
\partial_{\xi_n}\sigma_{-4}(D^{-3})](x_0)d\xi_n\sigma(\xi')dx' \nonumber\\
&=&i\int_{|\xi'|=1}\int^{+\infty}_{-\infty}{\rm trace} [\partial_{\xi_n}\pi^+_{\xi_n}\sigma_{-1}(D^{-1})\times
\sigma_{-4}(D^{-3})](x_0)d\xi_n\sigma(\xi')dx'.
\end{eqnarray}
In the normal coordinate, $g^{ij}(x_0)=\delta_i^j$ and $\partial_{x_j}(g^{\alpha\beta})(x_0)=0,$ {\rm if
}$j<n;~=h'(0)\delta^\alpha_\beta,~{\rm if }~j=n.$ So by Lemma A.2 in \cite{Wa3}, we have $\Gamma^n(x_0)=\frac{5}{2}h'(0)$ and
$\Gamma^k(x_0)=0$ for $k<n$. By the definition of $\delta^k$ and
Lemma 2.3 in \cite{Wa3}, we have $\delta^n(x_0)=0$ and $\delta^k=\frac{1}{4}h'(0)c(\widetilde{e_k})c(\widetilde{e_n})$ for $k<n$.
Then by (3.15), we obtain
\begin{eqnarray}
\sigma_{-4}(D^{-3})&=&\frac{c(\xi)\sigma_{2}(D^{3})c(\xi)}{|\xi|^8}-\frac{c(\xi)}{|\xi|^4}\sum_j\partial_{\xi_j}[c(\xi)|\xi|^2]
D_{x_j}[\frac{\sqrt{-1}c(\xi)}{|\xi|^4}]\nonumber\\
&=&\frac{1}{|\xi|^8}c(\xi)\Big(\frac{1}{2}h'(0)c(\xi)\sum_{k<n}\xi_k
c(\widetilde{e_k})c(\widetilde{e_n})-\frac{5}{2}h'(0)\xi_nc(\xi)-\frac{5}{4}h'(0)|\xi|^2c(dx_n)\Big)c(\xi)\nonumber\\
&&+\frac{\sqrt{-1}c(\xi)}{|\xi|^8}\Big(|\xi|^4c(dx_n)\partial_{x_n}c(\xi')-2h'(0)c(dx_n)c(\xi)
   +2\xi_{n}c(\xi)\partial_{x_n}c(\xi')+4\xi_{n}h'(0)\Big)\nonumber\\
 &=& \frac{1}{|\xi|^8}c(\xi)\Big(\frac{1}{2}h'(0)c(\xi)c(\xi')c(dx_n)-\frac{5}{2}h'(0)\xi_nc(\xi)-\frac{5}{4}h'(0)|\xi|^2c(dx_n)\Big)c(\xi)\nonumber\\
&&+\frac{1}{|\xi|^8}\Big(i|\xi|^4c(\xi)c(dx_n)\partial_{x_n}c(\xi')-2ih'(0)c(\xi)c(dx_n)c(\xi)\nonumber\\
  && +2i\xi_{n}c(\xi)c(\xi)\partial_{x_n}c(\xi')+4i\xi_{n}h'(0)c(\xi)\Big)\nonumber\\
 &=&\frac{1}{(1+\xi_{n}^{2})^{4}}\Big( (\frac{11}{2}\xi_{n}(1+\xi_{n}^{2})+8i\xi_{n})h'(0)c(\xi') \nonumber\\
 && +\big(-2i+6i\xi_{n}^{2}-\frac{7}{4}(1+\xi_{n}^{2})+\frac{15}{4}\xi_{n}^{2}(1+\xi_{n}^{2})\big)h'(0)c(dx_n) \nonumber\\
  &&-3i\xi_{n}(1+\xi_{n}^{2})\partial_{x_n}c(\xi')+i(1+\xi_{n}^{2})c(\xi')c(dx_n)\partial_{x_n}c(\xi') \Big).
\end{eqnarray}
 By (3.27), (3.31) and (3.36), we obtain
\begin{eqnarray}
&&{\rm tr} [\partial_{\xi_n}\pi^+_{\xi_n}\sigma_{-1}(D^{-1})\times
\sigma_{-4}(D^{-3})] \nonumber\\
&=&{\rm tr} \Big[-\frac{c(\xi')+ic(dx_n)}{2(\xi_n-i)^2} \times
\frac{1}{(1+\xi_{n}^{2})^{4}}\Big( \big(\frac{11}{2}\xi_{n}(1+\xi_{n}^{2})+8i\xi_{n}\big)h'(0)c(\xi') \nonumber\\
 && +\big(-2i+6i\xi_{n}^{2}-\frac{7}{4}(1+\xi_{n}^{2})+\frac{15}{4}\xi_{n}^{2}(1+\xi_{n}^{2})\big)h'(0)c(dx_n) \nonumber\\
  &&-3i\xi_{n}(1+\xi_{n}^{2})\partial_{x_n}c(\xi')+i(1+\xi_{n}^{2})c(\xi')c(dx_n)\partial_{x_n}c(\xi') \Big)\Big]\nonumber\\
&=&h'(0)\frac{7+6i-(20-15i)\xi_n+(7-6i)\xi_n^{2}+15i\xi_n^{3}}{(\xi_n-i)^5(\xi_n+i)^4}.
\end{eqnarray}
Then
\begin{eqnarray}
{\bf case~ b)}&=&
 i h'(0)\int_{|\xi'|=1}\int^{+\infty}_{-\infty}\frac{7+6i-(20-15i)\xi_n+(7-6i)\xi_n^{2}
   +15i\xi_n^{3}}{(\xi_n-i)^5(\xi_n+i)^4}d\xi_n\sigma(\xi')dx' \nonumber\\
&=& i h'(0)\frac{2 \pi i}{4!}\Big[\frac{7+6i-(20-15i)\xi_n+(7-6i)\xi_n^{2}+15i\xi_n^{3}}{(\xi_n+i)^4}
     \Big]^{(4)}|_{\xi_n=i}\Omega_4dx'\nonumber\\
&=&(-\frac{25}{8}-\frac{35i}{16})\pi h'(0)\Omega_4dx'.
\end{eqnarray}

 {\bf  case c)}~$r=-2,~l=-3,~k=j=|\alpha|=0$\\

From (2.4) we have
 \begin{equation}
{\rm case~ c)}=-i\int_{|\xi'|=1}\int^{+\infty}_{-\infty}{\rm trace} [\pi^+_{\xi_n}\sigma_{-2}(D^{-1})\times
\partial_{\xi_n}\sigma_{-3}(D^{-3})](x_0)d\xi_n\sigma(\xi')dx'.
\end{equation}
By (2.2.34) , (2.2.37) and (2.2.40) in \cite{Wa4}, we have
\begin{eqnarray}
\pi^+_{\xi_n}\sigma_{-2}(D^{-1})(x_0)|_{|\xi'|=1}
&=&\pi^+_{\xi_n}\Big[\frac{c(\xi)p_0(x_0)c(\xi)+c(\xi)c(dx_n)\partial_{x_n}[c(\xi')](x_0)}{(1+\xi_n^2)^2}\Big] \nonumber\\
&&-h'(0)\pi^+_{\xi_n}\Big[\frac{c(\xi)c(dx_n)c(\xi)}{(1+\xi_n)^3}\Big]\nonumber\\
&:=& B_1-B_2.
\end{eqnarray}
where
\begin{eqnarray}
B_1&=&\frac{-1}{4(\xi_n-i)^2}[(2+i\xi_n)c(\xi')p_0c(\xi')+i\xi_nc(dx_n)p_0c(dx_n) \nonumber\\
&&+(2+i\xi_n)c(\xi')c(dx_n)\partial_{x_n}c(\xi')+ic(dx_n)p_0c(\xi')
+ic(\xi')p_0c(dx_n)-i\partial_{x_n}c(\xi')]\nonumber\\
&=&\frac{1}{4(\xi_n-i)^2}\Big[\frac{5}{2}h'(0)c(dx_n)-\frac{5i}{2}h'(0)c(\xi')
  -(2+i\xi_n)c(\xi')c(dx_n)\partial_{\xi_n}c(\xi')+i\partial_{\xi_n}c(\xi')\Big]  ;         \\
B_2&=&\frac{h'(0)}{2}\Big[\frac{c(dx_n)}{4i(\xi_n-i)}+\frac{c(dx_n)-ic(\xi')}{8(\xi_n-i)^2}
+\frac{3\xi_n-7i}{8(\xi_n-i)^3}[ic(\xi')-c(dx_n)]\Big].
\end{eqnarray}
By (3.14), we obtain
 \begin{equation}
\partial_{\xi_n}\sigma_{-3}(D^{-3})=\frac{-4i\xi_nc(\xi')+(i-3i\xi_n^{2})c(dx_n)}{(1+\xi_n^{2})^3}.
\end{equation}
By (3.27), (3.42) and (3.43), we have
\begin{eqnarray}
&&{\rm tr }[B_2\times\partial_{\xi_n}\sigma_{-3}(D^{-3})(x_0)]|_{|\xi'|=1}\nonumber\\
&=&{\rm tr }\Big\{ \frac{h'(0)}{2}\Big[\frac{c(dx_n)}{4i(\xi_n-i)}+\frac{c(dx_n)-ic(\xi')}{8(\xi_n-i)^2}
+\frac{3\xi_n-7i}{8(\xi_n-i)^3}[ic(\xi')-c(dx_n)]\Big]\nonumber\\
&&\times \frac{-4i\xi_nc(\xi')+(i-3i\xi_n^{2})c(dx_n)}{(1+\xi_n^{2})^3}\Big\} \nonumber\\
&=&h'(0)\frac{4i-11\xi_n-6i\xi_n^{2}+3\xi_n^{3}}{(\xi_n-i)^5(\xi_n+i)^3}.
\end{eqnarray}
Similarly, we have
\begin{eqnarray}
&&{\rm tr }[B_1\times\partial_{\xi_n}\sigma_{-3}(D^{-3})(x_0)]|_{|\xi'|=1}\nonumber\\
&=&{\rm tr }\Big\{ \frac{1}{4(\xi_n-i)^2}\Big[\frac{5}{2}h'(0)c(dx_n)-\frac{5i}{2}h'(0)c(\xi')
  -(2+i\xi_n)c(\xi')c(dx_n)\partial_{\xi_n}c(\xi')+i\partial_{\xi_n}c(\xi')\Big]\nonumber\\
&&\times \frac{-4i\xi_nc(\xi')+(i-3i\xi_n^{2})c(dx_n)}{(1+\xi_n^{2})^3}\Big\} \nonumber\\
&=&h'(0)\frac{3+12i\xi_n+3\xi_n^{2}}{(\xi_n-i)^4(\xi_n+i)^3}.
\end{eqnarray}
Combining (3.44) and (3.45), we obtain
\begin{eqnarray}
{\bf case~ c)}&=&
 -i h'(0)\int_{|\xi'|=1}\int^{+\infty}_{-\infty}
 \frac{-7i+26\xi_n+15i\xi_n^{2}}{(\xi_n-i)^5(\xi_n+i)^3}d\xi_n\sigma(\xi')dx' \nonumber\\
&=&-i h'(0)\frac{2 \pi i}{4!}\Big[\frac{-7i+26\xi_n+15i\xi_n^{2}}{(\xi_n+i)^3}
     \Big]^{(5)}|_{\xi_n=i}\Omega_4dx'\nonumber\\
&=&\frac{55}{16}\pi h'(0)\Omega_4dx'.
\end{eqnarray}
Since $\Phi$ is the sum of the cases a), b) and c), so
  \begin{equation}
\Phi=(\frac{15}{16}-\frac{35i}{16})\pi h'(0)\Omega_4dx'.
\end{equation}

Now recall the Einstein-Hilbert action for manifolds with boundary \cite{Wa3}\cite{Wa4}\cite{Y},
 \begin{equation}
I_{\rm Gr}=\frac{1}{16\pi}\int_Ms{\rm dvol}_M+2\int_{\partial M}K{\rm dvol}_{\partial_M}:=I_{\rm {Gr,i}}+I_{\rm {Gr,b}},
\end{equation}
 where
 \begin{equation}
 K=\sum_{1\leq i,j\leq {n-1}}K_{i,j}g_{\partial M}^{i,j};~~K_{i,j}=-\Gamma^n_{i,j},
\end{equation}
 and $K_{i,j}$ is the second fundamental form, or extrinsic curvature. Taking the metric in Section 2, then by Lemma
A.2 \cite{Wa3}, for $n=6$, then
 \begin{equation}
K(x_0)=-\frac{5}{2}h'(0);~
 I_{\rm {Gr,b}}=-5h'(0){\rm Vol}_{\partial M}.
\end{equation}
 Then we obtain

\begin{thm}
Let  $M$ be a $6$-dimensional
compact spin manifold with the boundary $\partial M$ and the metric
$g^M$ as above and $D$ be the Dirac operator on $\widehat{M}$, then
\begin{equation}
{\rm Vol}^{(1,3)}_6=\widetilde{{\rm Wres}}[\pi^+D^{-1}\circ\pi^+D^{-3}]=-\frac{5\Omega_5}{3}\int_Ms{\rm dvol}_M
 +(\frac{7i}{8}-\frac{3}{8})\pi \Omega_4\int_{\partial M}Kd{\rm Vol}_{\partial M}.
\end{equation}
\end{thm}

\begin{rem}
In \cite{Wa3} \cite{Wa4}, Wang computed
 $\widetilde{{\rm Wres}}[\pi^+D^{-1}\circ\pi^+D^{-1}]$ and $\widetilde{{\rm Wres}}[\pi^+D^{-2}\circ\pi^+D^{-2}]$.
 In that cases, the boundary terms vanished, where the two operators are symmetric.  Theorem 3.4 states the
 boundary terms is non-zero when we compute $\widetilde{{\rm Wres}}[\pi^+D^{-1}\circ\pi^+D^{-3}]$. The reason
 is that $D^{-1}$ and $D^{-3}$ are not  symmetric.
\end{rem}

Let
 \begin{equation}
\widetilde{{\rm Wres}}[\pi^+D^{-1}\circ\pi^+D^{-3}]
=\widetilde{{\rm Wres}}_{i}[\pi^+D^{-1}\circ\pi^+D^{-3}]+\widetilde{{\rm Wres}}_{b}[\pi^+D^{-1}\circ\pi^+D^{-3}],
\end{equation}
where
 \begin{equation}
\widetilde{{\rm Wres}}_{i}[\pi^+D^{-1}\circ\pi^+D^{-3}]
=\int_M\int_{|\xi|=1}{\rm trace}_{S(TM)}[\sigma_{-6}(D^{-1-3})]\sigma(\xi)dx
\end{equation}
and
\begin{eqnarray}
&&\widetilde{{\rm Wres}}_{b}[\pi^+D^{-1}\circ\pi^+D^{-3}]\nonumber\\
&=& \int_{\partial M}\int_{|\xi'|=1}\int^{+\infty}_{-\infty}\sum^{\infty}_{j, k=0}\sum\frac{(-i)^{|\alpha|+j+k+1}}{\alpha!(j+k+1)!}
\times {\rm trace}_{S(TM)}[\partial^j_{x_n}\partial^\alpha_{\xi'}\partial^k_{\xi_n}\sigma^+_{r}(D^{-1})(x',0,\xi',\xi_n)
\nonumber\\
&&\times\partial^\alpha_{x'}\partial^{j+1}_{\xi_n}\partial^k_{x_n}\sigma_{l}(D^{-3})(x',0,\xi',\xi_n)]d\xi_n\sigma(\xi')dx'
\end{eqnarray}
denote the interior term and boundary term of $\widetilde{{\rm Wres}}[\pi^+D^{-1}\circ\pi^+D^{-3}]$.

Combining (3.48), (3.51) and (3.52), we obtain
\begin{cor}
Let  $M$ be a $6$-dimensional
compact spin manifold with the boundary $\partial M$ and the metric
$g^M$ as above and $D$ be the Dirac operator on $\widehat{M}$, then
\begin{eqnarray}
&& I_{\rm {Gr,i}}=\frac{-3}{80 \pi\Omega_5 }\widetilde{{\rm Wres}}_{i}[\pi^+D^{-1}\circ\pi^+D^{-3}]; \nonumber\\
&&I_{\rm {Gr,b}}=\frac{16}{(7i-3) \pi\Omega_4 }\widetilde{{\rm Wres}}_{b}[\pi^+D^{-1}\circ\pi^+D^{-3}].
\end{eqnarray}
\end{cor}

Nextly, for $5$-dimensional spin manifolds with boundary, we compute ${\rm Vol}^{(1,3)}_5$. By Proposition 2.2 (2) in \cite{Wa4}, we have
 \begin{equation}
\widetilde{{\rm Wres}}[\pi^+D^{-1}\circ\pi^+D^{-3}]=\int_{\partial M}\Phi.
\end{equation}
 By (2.4), when $n=5$, we have  $ r-k-|\alpha|+l-j-1=-5,~~r\leq -1,l\leq-3$, so we get $r=-1,~ l=-3,~k=|\alpha|=j=0,$ then
  \begin{equation}
\Phi=\int_{|\xi'|=1}\int^{+\infty}_{-\infty} {\rm trace}_{S(TM)}
[ \sigma^+_{-1}(D^{-1})(x',0,\xi',\xi_n)\times \partial_{\xi_n}\sigma_{-3}
(D^{-3})(x',0,\xi',\xi_n)]d\xi_n\sigma(\xi')dx'.
\end{equation}
By (2.2.44) in \cite{Wa3}, we have
  \begin{equation}
\pi^+_{\xi_n}\sigma_{-1} (D^{-1})(x_0)|_{|\xi'|=1}=\frac{c(\xi')+ic(dx_n)}{2(\xi_n-i)}.
\end{equation}
From (3.14) we obtain
 \begin{equation}
\partial_{\xi_n}\sigma_{-3}(D^{-3})=\frac{-4 i \xi_n c(\xi')+(i-3i\xi_n^{2})c(dx_n)}{(1+\xi_n^{2})^3}.
\end{equation}
Since $n=5$, ${\rm tr}(id)={\rm dim}(S(TM))=4$.
By the relation of the Clifford action and ${\rm tr}{AB}={\rm tr }{BA}$, then we have the equalities:
 \begin{equation}
{\rm tr}[c(\xi')c(dx_n)]=0;~~{\rm tr}[c(dx_n)^2]=-4;~~{\rm tr}[c(\xi')^2](x_0)|_{|\xi'|=1}=-4.
\end{equation}
Hence from (3.58), (3.59), and (3.60), we have
 \begin{eqnarray}
&&{\rm tr }\Big[\Big(\frac{c(\xi')+ic(dx_n)}{2(\xi_n-i)}\Big)\times
\Big(\frac{-4 i \xi_n c(\xi')+(i-3i\xi_n^{2})c(dx_n)}{(1+\xi_n^{2})^3}\Big)\Big]\nonumber\\
&=&\frac{2i-6\xi_n}{(\xi_n-i)^3(\xi_n+i)^3}.
\end{eqnarray}
Then
\begin{eqnarray}
\Phi&=& \int_{|\xi'|=1}\int^{+\infty}_{-\infty} \frac{2i-6\xi_n}{(\xi_n-i)^3(\xi_n+i)^3} d\xi_n\sigma(\xi')dx' \nonumber\\
&=& \frac{2 \pi i}{2!}\Big[\frac{2i-6\xi_n}{(\xi_n+i)^3} \Big]^{(2)}|_{\xi_n=i}\Omega_3dx'\nonumber\\
&=&\frac{3\pi i}{4} \Omega_3dx'.
\end{eqnarray}
By ${\rm Vol}^{(1,3)}_5=\frac{3\pi i}{4}\Omega_3{\rm Vol}_{\partial M}$ and
$I_{\rm {Gr,b}}=-4h'(0){\rm Vol}_{\partial M}$, we have

\begin{thm}
Let  $M$ be a $5$-dimensional
compact spin manifold with the boundary $\partial M$ and the metric
$g^M$ as in Section 2 and $D$ be the Dirac operator on $\widehat{M}$, then
\begin{eqnarray}
&& {\rm Vol}^{(1,3)}_5=\widetilde{{\rm Wres}}[\pi^+D^{-1}\circ\pi^+D^{-3}]=\frac{3\pi i}{4}\Omega_3{\rm Vol}_{\partial M}, \\
&& I_{\rm {Gr,b}}=\frac{16ih'(0)}{3\pi\Omega_3}\widetilde{{\rm Wres}}[\pi^+D^{-1}\circ\pi^+D^{-3}],
\end{eqnarray}
where ${\rm Vol}_{\partial M}$ denotes the canonical volume of ${\partial M}$.
\end{thm}

\section{A Kastler-Kalau-Walze type theorem for perturbations of Dirac operators}

 In \cite{CM}, Connes and Moscovici defined a twisted spectral triple and considered the perturbations of Dirac operator $e^{h}De^{h}$.
  In this section, for perturbations of Dirac operators, we compute the lower dimensional volume ${\rm Vol}^{(2,2)}_6$
 for $6$-dimensional spin manifolds with boundary and get a Kastler-Kalau-Walze type theorem in this case.

 \subsection{$\widetilde{{\rm Wres}}[\pi^+(fD^{-2})\circ\pi^+D^{-2}]$ for the Dirac operators}

Let $M$ be a $6$-dimensional compact spin manifold with the boundary $\partial M$ and the metric $g^M$ as Section 2 and
$D$ be the Dirac operator on $\widehat{M}$, we will compute $\widetilde{{\rm Wres}}[\pi^+(fD^{-2})\circ\pi^+D^{-2}]$
 for a smooth function $f(x)$.
 By the Kastler-Kalau-Walze type theorem in \cite{KW}, we get the lower dimensional volume for $6$-dimensional spin manifolds
 without boundary.  An application of Theorem 1 in \cite{Wa4} shows that

\begin{lem}
Let $M$ be a $6$-dimensional compact spin manifold without boundary, then
\begin{equation}
Wres[fD^{-4}]=-\frac{5\Omega_5}{3}\int_M fs{\rm dvol}_M.
\end{equation}
\end{lem}
Therefore, we only need to compute $\int_{\partial M}\Phi$. By Lemma 1 in \cite{Wa4}, we have
\begin{lem}
 For $6$-dimensional compact spin manifold with the boundary $\partial M$ and the metric $g^M$ as above, then
 \begin{eqnarray}
&&\sigma_{-2}(D^{-2})=|\xi|^{-2}; \\
&&\sigma_{-2}(fD^{-2})=f|\xi|^{-2}; \\
&&\sigma_{-3}(D^{-2})=-\sqrt{-1}|\xi|^{-4}\xi_k(\Gamma^k-2\delta^k)-\sqrt{-1}|\xi|^{-6}2\xi^j\xi_\alpha\xi_\beta
\partial_jg^{\alpha\beta}.
\end{eqnarray}
 \end{lem}
Since $\Phi$ is a global form on $\partial M$, so for any fixed point $x_0\in\partial M$, we can choose the normal coordinates
$U$ of $x_0$ in $\partial M$ (not in $M$) and compute $\Phi(x_0)$ in the coordinates $\widetilde{U}=U\times [0,1)\subset M$ and the
metric $\frac{1}{h(x_n)}g^{\partial M}+dx_n^2.$ For details, see Section 2.2.2\cite{Wa3}.
Now we can compute $\Phi$ (see formula (2.4) for the definition of $\Phi$), since the sum is taken over $
-r-l+k+j+|\alpha|=5,~~r,l\leq-2,$ then we have the following five cases:

{\bf case a)~I)}~$r=-2,~l=-2~k=j=0,~|\alpha|=1$

From (2.4) we have
 \begin{equation}
{\rm case~a)~I)}=-f\int_{|\xi'|=1}\int^{+\infty}_{-\infty}\sum_{|\alpha|=1}
{\rm trace}[\partial^\alpha_{\xi'}\pi^+_{\xi_n}\sigma_{-2}(D^{-2})\times
\partial^\alpha_{x'}\partial_{\xi_n}\sigma_{-2}(D^{-2})](x_0)d\xi_n\sigma(\xi')dx'.
\end{equation}
By Lemma 3.2, for $i<n$, then
\begin{eqnarray}
\partial_{x_i}\sigma_{-2}(D^{-2})(x_0)=\partial_{x_i}{(|\xi|^{-2})}(x_0)=
-\frac{\partial_{x_i}(|\xi|^2)(x_0)}{|\xi|^4}=0.
\end{eqnarray}
Then case a) I) vanishes.

 {\bf case a)~II)}~$r=-2,~l=-2~k=|\alpha|=0,~j=1$

From (2.4) we have
\begin{eqnarray}
{\rm case~a)~II)}&=&-\frac{1}{2}\int_{|\xi'|=1}\int^{+\infty}_{-\infty} {\rm
trace} [\partial_{x_n}\pi^+_{\xi_n}\sigma_{-2}(fD^{-2})\times\partial_{\xi_n}^2\sigma_{-2}(D^{-2})](x_0)d\xi_n\sigma(\xi')dx' \nonumber\\
&=&-\frac{1}{2}f\int_{|\xi'|=1}\int^{+\infty}_{-\infty} {\rm
trace} [\partial_{x_n}\pi^+_{\xi_n}\sigma_{-2}(D^{-2})\times\partial_{\xi_n}^2\sigma_{-2}(D^{-2})](x_0)d\xi_n\sigma(\xi')dx'\nonumber\\
&&-\frac{1}{2}\partial_{x_{n}}(f)\int_{|\xi'|=1}\int^{+\infty}_{-\infty} {\rm
trace} [\pi^+_{\xi_n}\sigma_{-2}(D^{-2})\times\partial_{\xi_n}^2\sigma_{-2}(D^{-2})](x_0)d\xi_n\sigma(\xi')dx'.
\end{eqnarray}
By {\bf case a)~II)} in \cite{Wa4}, we have
\begin{equation}
-\frac{1}{2}f\int_{|\xi'|=1}\int^{+\infty}_{-\infty} {\rm
trace} [\partial_{x_n}\pi^+_{\xi_n}\sigma_{-2}(D^{-2})\times\partial_{\xi_n}^2\sigma_{-2}(D^{-2})](x_0)d\xi_n\sigma(\xi')dx'
=-\frac{5}{8}\pi f h'(0)\Omega_4dx',
\end{equation}
where $\Omega_4$ is the canonical volume of $S^4$.

On the other hand, by (14) in \cite{Wa4}, we have
\begin{equation}
\partial^2_{\xi_n}(\sigma_{-2}(D^{-2}))(x_0)=\partial^2_{\xi_n}(|\xi|^{-2})(x_0)=\frac{-2+6\xi_n^2}{(1+\xi_n^2)^3}.
\end{equation}
By (4.2) and the Cauchy integral formula , then
\begin{eqnarray}
\pi^+_{\xi_n}\sigma_{-2}(D^{-2})(x_0)|_{|\xi'|=1}&=&
 \frac{1}{2\pi i}\lim_{u\rightarrow0^-}\int_{\Gamma^+}\frac{\frac{1}{(\eta_n+i)(\xi_n+iu-\eta_n)}}
   {(\eta_n-i)}d\eta_n  \nonumber\\
&=& \frac{-i}{2(\xi_n-i)}.
\end{eqnarray}
Combining (4.9) and (4.10), we have
\begin{equation}
{\rm tr }\Big[\frac{-i}{2(\xi_n-i)} \times \frac{-2+6\xi_n^2}{(1+\xi_n^2)^3}\Big]
=\frac{i-3i\xi_n^{2}}{(\xi_n-i)^4(\xi_n+i)^3}{\rm tr }[id]=\frac{8i-24i\xi_n^{2}}{(\xi_n-i)^4(\xi_n+i)^3},
\end{equation}
where $n=6$, ${\rm tr}_{S(TM)}[{\rm id}]={\rm dim}(\wedge^*(3))=8$.

From (4.7) and (4.11), we have
\begin{eqnarray}
&&-\frac{1}{2}\partial_{x_{n}}f\int_{|\xi'|=1}\int^{+\infty}_{-\infty} {\rm
trace} [\pi^+_{\xi_n}\sigma_{-2}(D^{-2})\times\partial_{\xi_n}^2\sigma_{-2}(D^{-2})](x_0)d\xi_n\sigma(\xi')dx'\nonumber\\
&=& -\frac{1}{2}\partial_{x_{n}}(f)\int_{|\xi'|=1}\int^{+\infty}_{-\infty} \frac{8i-24i\xi_n^{2}}{(\xi_n-i)^4(\xi_n+i)^3}
 (x_0)d\xi_n\sigma(\xi')dx' \nonumber\\
&=& -\frac{1}{2}\partial_{x_{n}}(f)\frac{2 \pi i}{3!}\Big[\frac{8i-24i\xi_n^{2}}{(\xi_n+i)^3} \Big]^{(3)}|_{\xi_n=i} \Omega_3dx'\nonumber\\
&=&3\pi i \partial_{x_{n}}(f) \Omega_4dx'.
\end{eqnarray}
Hence in this case,
\begin{equation}
{\rm case~a)~II)}=-\frac{5}{8}\pi f h'(0)\Omega_4dx'+3\pi i \partial_{x_{n}}(f) \Omega_4dx'.
\end{equation}

{\bf case a)~III)}~$r=-2,~l=-2~j=|\alpha|=0,~k=1$

From (2.4) we have
\begin{eqnarray}
{\rm case~ a)~III)}&=&-\frac{1}{2}\int_{|\xi'|=1}\int^{+\infty}_{-\infty}
{\rm trace} [\partial_{\xi_n}\pi^+_{\xi_n}\sigma_{-2}(fD^{-2})\times
\partial_{\xi_n}\partial_{x_n}\sigma_{-2}(D^{-2})](x_0)d\xi_n\sigma(\xi')dx' \nonumber\\
&=&-\frac{1}{2}f\int_{|\xi'|=1}\int^{+\infty}_{-\infty}
{\rm trace} [\partial_{\xi_n}\pi^+_{\xi_n}\sigma_{-2}(D^{-2})\times
\partial_{\xi_n}\partial_{x_n}\sigma_{-2}(D^{-2})](x_0)d\xi_n\sigma(\xi')dx'.
\end{eqnarray}
By {\bf case a)~III)} in \cite{Wa4}, we have
\begin{equation}
-\frac{1}{2}\int_{|\xi'|=1}\int^{+\infty}_{-\infty}{\rm trace} [\partial_{\xi_n}\pi^+_{\xi_n}\sigma_{-2}(D^{-2})\times
\partial_{\xi_n}\partial_{x_n}\sigma_{-2}(D^{-2})](x_0)d\xi_n\sigma(\xi')dx'
=\frac{5}{8}\pi h'(0)\Omega_4dx'.
\end{equation}
Combining (4.14) and (4.15), we obtain
\begin{equation}
{\rm case~ a)~III)}=\frac{5}{8}\pi f h'(0)\Omega_4dx'.
\end{equation}

Since {\bf case b)}, {\bf case c)} has the same expression with the case of {\bf case b)}, {\bf case c)} in \cite{Wa4}
multiplied by a function $f$ ,
 so we can use the same way to compute the two terms. An application of (21) and (24) in \cite{Wa4} shows that the
 sum of  {\bf case b)} and {\bf case c)} is zero. Hence we conclude that,  the sum of {\bf case~
a)}, {\bf case b)} and {\bf case c)} is
\begin{equation}
\Phi=3\pi i \partial_{x_{n}}(f) \Omega_4dx'.
\end{equation}

In summary, we have proved
\begin{thm}
Let $M$ be a $6$-dimensional compact spin manifold with the boundary $\partial M$ and the metric
$g^M$ as above and $D$ be the Dirac operator on $\widehat{M}$, then
\begin{equation}
\widetilde{{\rm Wres}}[\pi^+(fD^{-2})\circ\pi^+D^{-2}]=-\frac{5\Omega_5}{3}\int_Mfs{\rm dvol}_M
  +3\pi i  \Omega_4\int_{\partial M}\partial_{x_{n}}(f)|_{x_{n}=0} {\rm dvol}_{\partial M}.
\end{equation}
\end{thm}

\subsection{$\widetilde{{\rm Wres}}[\pi^+(D^{2}+f)^{-1}\circ\pi^+D^{-2}]$ for the Dirac operators}
In this subsection, let $M$ be a $6$-dimensional compact spin manifold with the boundary $\partial M$ and the metric $g^M$ as Section 2 and
$D$ be the Dirac operator on $\widehat{M}$,  we will compute $\widetilde{{\rm Wres}}[(D^{2}+f)^{-1}D^{-2}]$ for  smooth function $f(x)$.

Firstly, we compute the symbol $\sigma((D^{2}+f)^{-1})$ of $(D^{2}+f)^{-1}$.
By the definition of the Dirac operator $D$  in section 3 and Lemma 1 in \cite{Wa4}, we have
\begin{equation}
\sigma_{-2}((D^{2}+f)^{-1})=\sigma_{-2}(D^{-2});~~~\sigma_{-3}((D^{2}+f)^{-1})=\sigma_{-3}(D^{-2}).
\end{equation}
Let $\tilde{D}=D^{4}+fD^{2}$ be the fourth order operators. Write
 \begin{equation}
D_x^{\alpha}=(-\sqrt{-1})^{|\alpha|}\partial_x^{\alpha};
~\sigma(\tilde{D})=\tilde{p}_4+\tilde{p}_3+\tilde{p}_2+\tilde{p}_1+\tilde{p}_0;
~\sigma((\tilde{D})^{-1})=\sum^{\infty}_{j=4}\tilde{q}_{-j}.
\end{equation}
Since the smooth function $f(x)\in \tilde{p}_0$, then we have
\begin{equation}
\sigma_{4}(\tilde{D}^{4})=\sigma_{4}(D^{4});~~\sigma_{3}(\tilde{D}^{4})=\sigma_{3}(D^{4});
~~\sigma_{2}(\tilde{D}^{4})=\sigma_{2}(D^{4})+f|\xi|^2.
\end{equation}
By the composition formula of psudodifferential operators, we have
 \begin{eqnarray}
1=\sigma(\tilde{D} \circ \tilde{D}^{-1})
&=&\sum_{\alpha}\frac{1}{\alpha!}\partial^{\alpha}_{\xi}[\sigma(\tilde{D})]D^{\alpha}_{x}[\sigma((\tilde{D})^{-1})] \nonumber\\
&=&(\tilde{p}_4+\tilde{p}_3+\tilde{p}_2+\tilde{p}_1+\tilde{p}_0)(\tilde{q}_{-4}+\tilde{q}_{-5}+\tilde{q}_{-6}+\cdots) \nonumber\\
&&+\sum_j(\partial_{\xi_j}\tilde{p}_4+\partial_{\xi_j}\tilde{p}_3+\partial_{\xi_j}\tilde{p}_2
+\partial_{\xi_j}\tilde{p}_1+\partial_{\xi_j}\tilde{p}_0) \nonumber\\
&& \times \Big(D_{x_j}\tilde{q}_{-4}+D_{x_j}\tilde{q}_{-5}+D_{x_j}\tilde{q}_{-6}+\cdots\Big) \nonumber\\
&&+\sum_{i,j}\Big(\partial_{\xi_i}\partial_{\xi_j}(\tilde{p}_4+\tilde{p}_3+\tilde{p}_2+\tilde{p}_1+\tilde{p}_0) \Big)\nonumber\\
&& \times \Big( D_{x_i}D_{x_j}\tilde{q}_{-4}+D_{x_i}D_{x_j}\tilde{q}_{-5}+D_{x_i}D_{x_j}\tilde{q}_{-6}+\cdots\Big) \nonumber\\
\nonumber\\
&=&\tilde{p}_4\tilde{q}_{-4}+\Big(\tilde{p}_4\tilde{q}_{-5}+\tilde{p}_3\tilde{q}_{-4}
+\sum_j\partial_{\xi_j}\tilde{p}_4D_{x_j}\tilde{q}_{-4}\Big)+\cdots.
\end{eqnarray}
Then we obtain
\begin{eqnarray}
\tilde{q}_{-4}&=&\tilde{p}_4^{-1}; \\
\tilde{q}_{-5}&=&-\tilde{p}_4^{-1}[\tilde{p}_3\tilde{q}_{-4}+\sum_j\partial_{\xi_j}\tilde{p}_4D_{x_j}(\tilde{q}_{-4})]; \\
0&=&\tilde{p}_4\tilde{q}_{-6}+\tilde{p}_3\tilde{q}_{-5}+\tilde{p}_2\tilde{q}_{-4}
  +\sum_j\partial_{\xi_j}\tilde{p}_3D_{x_j}\tilde{q}_{-4}+\sum_j\partial_{\xi_j}\tilde{p}_4D_{x_j}\tilde{q}_{-5}\nonumber\\
&&+\sum_{i,j}\partial_{\xi_i}\partial_{\xi_j}\tilde{p}_4D_{x_i}D_{x_j}\tilde{q}_{-4} .
\end{eqnarray}
Hence from (4.21)-(4.25) and the recursion formula, we have

\begin{lem}
For a $6$-dimensional compact spin manifold with the boundary $\partial M$, then
\begin{equation}
q_{-6}=-f|\xi|^2+\sigma_{-6}(D^{-4}).
\end{equation}
\end{lem}
By Lemma 4.4 and \cite{Ka}, \cite{KW}, we obtain

\begin{thm}
let $M$ be a $6$-dimensional compact spin manifold without boundary and $f_{i}(i=1,2,3)$ be the smooth function, then
\begin{equation}
{\rm Wres}[(D^{4}+f_{1}D^{2}+f_{2}D+f_{3})^{-1}]=-\Omega_5\int_M(\frac{5}{3}s+8f_{1}){\rm dvol}_M.
\end{equation}
\end{thm}
\begin{rem}
When $f_{i}=0~(i=1,2,3)$, we get the classical Kastler-Kalau-Walze type theorem.
\end{rem}
On the other hand, it is straightforward to see from (4.19) we can directly get the same results with
Theorem 1 in \cite{Wa4}. Hence we conclude that $\int_{\partial M}\Phi=0$. In summary, we have proved

\begin{thm}
let $M$ be a $6$-dimensional compact spin manifold with the boundary $\partial M$ and the metric $g^M$ as Section 2 and
$D$ be the Dirac operator on $\widehat{M}$, then
\begin{equation}
\widetilde{{\rm Wres}}[\pi^+(D^{2}+f)^{-1}\circ\pi^+D^{-2}]=-\Omega_5\int_M(\frac{5}{3}s+8f){\rm dvol}_M.
\end{equation}
\end{thm}

\section*{ Acknowledgements}
This work was supported by Fok Ying Tong Education Foundation under Grant No. 121003 and NSFC. 11271062. The author also thank the referee
for his (or her) careful reading and helpful comments.

\end{document}